\documentclass{article}
\usepackage{amsmath,amsthm,amsfonts,graphicx}
\usepackage{multirow}
\usepackage{slashbox}
\usepackage{multirow}
\def\smallddots{\mathinner{\raise7pt\hbox{.}\raise4pt\hbox{.}\raise1pt\hbox{.}}}
\def\smallsdots{\mathinner{\raise1pt\hbox{.}\raise4pt\hbox{.}\raise7pt\hbox{.}}}

\numberwithin{equation}{section}
\numberwithin{table}{section}
\newtheorem{theorem}{Theorem}

\newtheorem{corollary}{Corollary}

\newtheorem{remark}{Remark}

\begin{document}

{\Large \bf Better Late Than Never: Filling a Void in the History

\centerline{of Fast Matrix Multiplication and Tensor Decompositions}}

\medskip

\medskip

\centerline{Victor Y. Pan}

\medskip

\medskip

\centerline{Department of Mathematics and Computer Science}

\centerline{Lehman College of the City University of New York}

\centerline{Bronx, NY 10468 USA}

\centerline{and}
\centerline{Ph.D. Programs in Mathematics  and Computer Science}
\centerline{The Graduate Center of the City University of New York}
\centerline{New York, NY 10036 USA}
\centerline{victor.pan@lehman.cuny.edu}
\centerline{http://comet.lehman.cuny.edu/vpan/}

 \date{}


\begin{abstract}
Multilinear and tensor decompositions are a popular tool
in linear and multilinear algebra and have a wide range 
of important applications to modern computing. Our paper 
of 1972 presented the first nontrivial application of such 
decompositions to fundamental matrix computations and was 
also a landmark in the history of the acceleration of matrix
multiplication. Published in 1972 in Russian, it has never 
been translated into English. It has been very rarely cited 
in the Western literature on matrix multiplication and never 
in the works on multilinear and tensor decompositions. This 
motivates us to present its translation into English, together 
with our brief comments on its impact on the two fields.
\end{abstract}


\paragraph{\bf 2000 Math. Subject Classification:}
15A69, 01A60, 15-03

\paragraph{\bf Key Words:}
Tensor decomposition,
Multilinear algebra,
Matrix multiplication,
Linear algebra

\bigskip

Multilinear or tensor decompositions have been introduced by Hitchcock
in 1927, but received scant attention until
a half of a dozen of papers appeared in 1963--70
in the psychometrics literature.
Our paper \cite{P72}
was the next major step
in this history.
It presented the earliest known application of nontrivial
multilinear and tensor decompositions to fundamental matrix computations.
By now such decompositions have become a popular tool
in linear and multilinear algebra and have a wide range of
important applications to modern computing
(see \cite{T03}, \cite{KB09},  \cite{OT10}, and the bibliography therein).
The paper  \cite{P72} was also a landmark in the history of the acceleration of matrix
multiplication, hereafter referred to as MM
(see our recent brief review of the history of that study in \cite{P14}).

Let us supply some specific comments.
Hereafter $MM(m,k,n)$ denotes the problem of
multiplying a matrix
 $A$ of size $m\times k$ by a matrix $B$ of size $k\times n$,
which are square matrices if $m=k=n$.
Let $D$ denote an auxiliary  $n\times m$ matrix
and let $l_s(A)$,
$l'_s(B)$, and $l''_s(D)$ denote
 linear forms
in the entries of the matrices $A$, $B$, and  $D$,
respectively.

According to Part 1 of \cite[Theorem 2]{P72},
we can perform $M(K,K,K)$ by using $cK^{\omega}$
arithmetic operations for any positive integer $K$,
 a constant $c$ independent of $K$, and $\omega=3\log_{mkn} r_{mkn}$
provided that we are given a bilinear algorithm of rank $r_{mkn}$
for $M(m,k,n)$, that is, a bilinear decomposition of the set 
of the $mn$ entries of the product $AB$ into the sum of $r_{mkn}$ 
bilinear products.
This reduces square $MM(K,K,K)$ for all $K$ 
to rectangular $M(m,k,n)$ for any fixed triple of $m$,
 $k$, and  $n$.

\cite[Theorem 2]{P72} reduces bilinear decomposition of rank $R$ for
$MM(m,k,n)$ 
 to decomposition of
the trilinear form trace$(ABD)$ of rank  $R$, that is,
its decomposition
into a sum $\sum_{s=1}^R l_s(A)l'_s(B)l''_s(D)$.
 That simple but basic result turns the problem 
of the acceleration of MM  into the search for
trilinear decompositions of trace$(ABD)$
of smaller rank  
where $A$, $B$, and $D$ denote 3 matrices
of fixed sizes.
The paragraph following that theorem of \cite{P72}
demonstrates a novel nontrivial
technique of trilinear aggregation for generating
such decompositions. The demonstration is by presenting
a decomposition  of
trace$(ABD)$ for $MM(n,n,n)$ generated by means of using this technique.
The decomposition has rank $0.5n^3+3n^2$ for any even $n$
(versus the straightforward decomposition of
rank $n^3$).
After its refinement in 1978 in  \cite{P78},
the technique of trilinear aggregation
enabled a decomposition of rank $\frac{1}{3}(n^3-4n)+6n^2$
for the trace$(ABD)$ for $MM(n,n,n)$.
For $n=70$ the rank turns into 143,640, and this implied
the
decrease of the record exponent of MM, $\log_2(7)<2.8074$, established
in 1969 in \cite{S69}, to $3\log_{70}143,640<2.7962$, that is, 
implied the acceleration of $MM(n,n,n)$
from $cn^{2.8704}$ for a constant $c$ to $O(n^{2.7962})$.

This was a landmark in the study of MM because since  1969,
when the MM exponent was decreased  from 3 to 2.8074 in \cite{S69},
all leading experts worldwide intensified their effort  in competition 
for decreasing it further, towards the information lower bound 2. 
 The exponent 2.8074, however,  has defied these intensive attacks 
 for almost a decade, from 1969 to 1978, until it
was beaten in  \cite{P78}, based on the application of 
trilinear aggregation. The new exponent
was decreased a number of times soon thereafter, and
trilinear aggregation
remained an indispensable ingredient of almost all
algorithms supporting this development
(see \cite[page 255]{CW90}).

The other results of the paper \cite{P72} were also of interest 
for the study of MM.
\cite[Theorem 1]{P72} has established some early upper and lower bounds
on the arithmetic complexity
of the general problem $MM(m,k,n)$ for any triple of $m$,
$k$, and $n$ and on its rank, that is, the minimal rank of its bilinear 
and trilinear decompositions.
The theorem has also
 estimated the arithmetic complexity of MM
and its rank in terms of one another
and has bounded
the rank and arithmetic complexity
 of the specific problems
 $MM(2,2,n)$, $MM(2,3,3)$, $MM(2,3,3)$,  $MM(2,3,4)$, and
$MM(2,4,4)$. Moreover part 5 of \cite[Theorem 1]{P72}
showed the duality of 6 algorithms for the 6 problems
$MM(m,k,n)$, $MM(k,n,m)$, $MM(n,m,k)$, $MM(k,m,n)$, $MM(n,k,m)$,
and $MM(m,n,k)$.
In \cite{HM73} and \cite{W80} the duality 
technique was applied to the
the design of efficient algorithms for MM, convolution,
integer multiplication, and the
FIR filters.
\cite[Theorem 3]{P72} provides a nontrivial explicit expression for all
bilinear algorithms of rank 7 for $MM(2,2,2)$ (cf. \cite{dG78}).

The paper \cite{P72} has been written by the author in his 
spare time while he was working in Economics to make his living. 
The paper was published in Russian, never translated into English so far,
 very rarely cited in the Western literature on MM
and never
in the works
on multilinear and tensor decompositions.
This motivates us to undertake its
 translation into English,
which we present below.

\bigskip

{\bf Acknowledgements:}
Our work has been supported by NSF Grant CCF--1116736 and
PSC CUNY Award 67699--00 45.



\clearpage

``USPEKHI MATEMATICHESKIKH NAUK"
(``RUSSIAN MATH. SURVEYS"),

\medskip

\centerline{
XXVII, 5(167), 249--250,
Sept.-- Oct. 1972}

\medskip
\medskip
\medskip

In the MOSCOW MATHEMATICAL SOCIETY

\medskip
\medskip

\medskip

\centerline{ON ALGORITHMS FOR MATRIX MULTIPLICATION AND INVERSION}

\medskip
\medskip
\medskip

\centerline{V. Y. Pan}

\medskip
\medskip
\medskip

Suppose we seek the values of $M$ rational functions
$R_l=R_l (x_1,\ldots,x_N)$, $l=1,\dots,M$, at $N$ points
 $x_1,\ldots,x_N$. Fix a sequence of arithmetic operations
 $p_i=R_i^{'}~o~ R_i^{"}$, $i=1,\ldots,Q$, such that
 every symbol $o$ stands for $+$, $-$, $\times$ or $\div$;
every $R_i^{'}$ as well as $R_i^{"}$
(for $i=1,\dots,Q$) is either a constant, or one of the $x_1,\ldots,x_N$, or  $p_j$ for $j<i$,
and $p_{s_l}\equiv R_l$  identically in  $x_1,\ldots,x_N$ for $1\le s_l\le Q$ and $l=1,\ldots,M$.
Call such a sequence a {\em computational program}
of {\em length} $Q$ for computing
the functions $R_1,\ldots,R_M$ (see \cite[page 104]{P1966}).
See some examples of these programs for $N=1$ in \cite{P1966}.

Consider such programs for computing a matrix product
$C_{mn}=A_{mk} B_{kn}$ and the inverse matrix $D_n=(A_n)^{-1}$.
In this case the inputs $x_h$ are given by the entries of the matrices
$A_{mk}$, $B_{kn}$, and $A_n$, and the outputs $R_l$ are
the entries of the matrices
 $C_{mn}$ and $D_n$. Let $Q_{mkn}$ and $Q_n^{-}$ denote the lengths of
the programs.
We naturally arrive at the two problems
of devising programs of the minimal lengths for  matrix multiplication and inversion.
Below we reduce both of them to a single problem
expressed through identity (\ref{eq1}) below
(see parts 1--3 of Theorem \ref{th1}).

\

Given four positive integers $m$, $k$, $n$, and $R$
and  three matrices $\{U,V,W\}=\{u_{ij}^{(s)},v_{gh}^{(s)},w_{lq}^{(s)}\}$ for
$g,j=0,1,\ldots,k-1$; $i,l=0,1,\ldots,m-1$; $h,q=0,1,\ldots,n-1$; $s=1,\ldots,R$,
assume the following identities in $a_{ij}$ and $b_{gh}$,

\begin{equation}\label{eq1}
c_{lq}=\Sigma_{g=0}^{k-1}a_{lg}b_{gq}\equiv
\Sigma_{s=1}^Rw_{lq}^{(s)} (\Sigma_{i=0}^{m-1}\Sigma_{j=0}^{k-1}u_{ij}^{(s)}  a_{ij})(\Sigma_{g=0}^{k-1}\Sigma_{h=0}^{n-1}v_{gh}^{(s)}  b_{gh})
\end{equation}

\centerline{for $l=0,\dots,m-1$ and $q=0,\dots,n-1$.}
\noindent Then we call (\ref{eq1}) {\em elementary program} for multiplication of matrices $A_{mk}$ and $B_{kn}$,
 call the numbers  $u_{ij}^{(s)},v_{gh}^{(s)},w_{lq}^{(s)}$ the elements of the program,
and call $R=R_{mkn}$ its length.
Let $q_{mkn}$ and $q_n^{-}$ denote the minimums of the lengths $Q_{mhn}$ and
$Q_n^{-}$,
respectively, over all computational programs and let  $r_{mkn}$
denote the minimum of the length $R_{mkn}$
over all
elementary programs.

\

Clearly, identity (\ref{eq1}) still holds where $c_{lq},a_{ij},b_{gh}$
are not numbers but, say, square matrices of
a fixed size. For their multiplication we can again apply elementary program (\ref{eq1}).
By repeating this recursively, we can reduce
general matrix multiplication to multiplication of matrices of small
sizes. More precisely, every decrease of the integer value
$p=\displaystyle\min_{m,k,n}p_{mkn}$, for $p_{mkn}=3 \log_{mkn}r_{mkn}$,
implies a decrease of the order of the length $q_{KKK}$ as $K\rightarrow \infty$
(see  part 1 of Theorem \ref{th1} below).
The construction in \cite{S1969} of an elementary program
supporting $R_{222}=7$ has immediately yielded programs of
length $Q_{KKK}\le C(2,2,2) K^{\log_{2}{7}} $  for all $K$.
Furthermore, the orders of magnitude of the values $r_{mkn}$  and   $q_{mkn}$  coincide
with one another
up to within the term  $mk+kn+mn$   (see below part 3 of Theorem \ref{th1}).
Let us state these facts formally and complement them with some estimates for  $r_{mkn}$
and  $q_{mkn}$.

\begin{theorem}\label{th1}
Fix four natural numbers $K$, $m$, $k$, and $n$. Then

    1)	(V. Strassen)~ $Q_{KKK}\le C(m,k,n)
K^{p_{mkn}}$;~~

$r_{222}\leq 7
	\rightarrow q_{KKK}\leq C(2,2,2) K^{\log_{2}{7}} .$
	
\
	
    2)~$\frac{1}{6} (q_{2n}^{-} -2q_{n}^{-}-3n^2 )\le q_{nnn} \le 4q_n^{-} +2q_{2n}^{-}+n^2+n.$

\

3)~$(m+n-1)k\le r_{mkn}\le C_{1}~q_{mkn}+C_2 (mk+kn+mn).$
	
\

	4)~$q_{mkn}\ge 2(m+n-1)k-m-n+1.$
	
\
	
	5)~$r_{mkn}=r_{mnk}=r_{nkm}=r_{nmk}= r_{knm}= r_{kmn}.$
	
\

	6)~$r_{22n}\ge 3n+2~({\rm for}~ n\ge 3);~ r_{222}\ge 7~({\rm see}~ [3]);
~15\le r_{233}\le 16$;

$r_{234}\ge 19; ~r_{333}\ge 18; ~ r_{244}\le 27.$
\end{theorem}

\noindent Here $C_{1}$ and $C_{2}$ are positive constants; $C(m,k,n)$ does not depend on $K$.

\begin{remark}\label{re1}
The statement in part 4 of the theorem can be made more precise:
every program that computes $c_{0h}=\Sigma_{g=0}^{k-1}a_{0g} b_{gh}; c_{l0}=\Sigma_{g=0}^{k-1}a_{lg} b_{g0}
 (l=0,\dots,m-1;h=0,\dots,n-1)$ for the input
 $a_{lg},b_{gh}$ ($l=0,\dots,m-1$; $h=0,\dots,n-1$; $g=0,\dots,k-1$),  must use at least
 $(m+n-1)k$ multiplications and divisions
and at least $(m+n-1)k-m-n-1$ additions or subtractions
(this is proved by methods that generalize the proof of \cite[Theorem 1.1]{P1966}).
 \end{remark}

 \begin{corollary}\label{co1}
 $p_{22n}\geq p_{222}$ (for $n=1,2,3,\dots);~ p_{233}> p_{222}.$
\end{corollary}
Next we show an alternative method (distinct from the one of \cite{S1969}) for fast multiplication of matrices.

\begin{theorem}\label{th2}
The entries of the three matrices $U$, $V$, and $W$ support an elementary program (\ref{eq1}) of  length $R$
 if and only if the following identity in $a_{ij},b_{gh},d_{ql}$ holds,
$$\Sigma_{s=1}^{R}(\Sigma_{i,j}u_{ij}^{s})
(\Sigma_{g,h}v_{gh}^{s} b_{gh})(\Sigma_{l,q} u_{lq}^{s} d_{ql})
\equiv \Sigma_{i,j,h} a_{ij} b_{jh} d_{hi}.$$
\end{theorem}

The following identity in $a_{ij},b_{gh},d_{ql}$
defines an elementary  program (\ref{eq1})
where $R_{nnn}=0.5n^3+3n^2$ and $n=2m$ (for $m=1,2,\dots$),

$$\Sigma_{i,j,h=0;~i+j+h~{\rm is~even}}^{n-1}(a_{ij}+a_{h+1,i+1})(b_{jh}+b_{i+1,j+1})(d_{hi}+d_{j+1,h+1})-$$

$$\Sigma_{i,h=0}^{n-1}a_{h+1,i+1}\Sigma_{j=0;~j:~i+j+h~{\rm is~even}}^{n-1}(b_{jh}+b_{i+1,j+1})  d_{hi}-$$

$$\Sigma_{i,h=0}^{n-1}a_{ij} b_{i+1,j+1}            \Sigma_{h=0;~h:~i+j+h~{\rm is~even}}^{n-1}(d_{hi}+d_{j+1,h+1})-$$

 $$\Sigma_{h,j=0}^{n-1}\Sigma_{i=0;~i:~i+j+h~{\rm is~even}}^{n-1}(a_{ij}+a_{h+1,i+1} ) b_{jh} d_{j+1,h+1}\equiv
\Sigma_{i,j,h}a_{ij} b_{jh} d_{hi}.$$
Here $f_{ln}=f_{l0}$ and $f_{nl}=f_{0l}$ for $l=0,1,\dots,n$, while $f$ stands for $a;b;d.$
It follows that $P_{34,34,34}\approx 2.8495$, that is,
this method multiplies matrices by using less than $O(K^{3})$ operations.

Let us establish a uniqueness property of the algorithm of the paper
\cite{S1969} for computing $A_{22} B_{22}$.
Two triples $\{U,V,W\}$ and $\{\overline {U},\overline{V},\overline{W} \}$
defining an elementary scheme (\ref{eq1})
are said to be equivalent
to one another if
$$\overline{u}_{ij}^{s}=
\Sigma_{v,\kappa}\sigma_{iv} \nabla_{j\kappa} u_{v\kappa}^{t(s)},~\overline{v}_{gh}^{(s)}=
\Sigma_{v,\kappa}\lambda_{vg} \mu_{h\kappa} v_{v\kappa}^{t(s)},~{\rm and}~\overline{w}_{lq}^{s}=
\Sigma_{v,\kappa}\gamma_{vl}\beta_{\kappa q}w_{v\kappa}^{t(s)},$$
where the matrices in the three pairs $(\sigma_{iv})$ and $(\gamma_{vl} ),(\nabla_{j\kappa})$ and $(\lambda_{vg} )$,
and $(\mu_{h\kappa})$
and $(\beta_{\kappa q} )$ are the inverses of one another; $1\leq t(s)  \leq R$;
$t_{s_{1}} \neq t_{s_{2}}$ for $s_{1} \neq {s_{2}}$, and all $t_{s}$ are integers.

\begin{theorem}\label{th3}
An elementary program for computing the matrix product
$A_{22} B_{22}$ has a length $R_{222} \leq 7$
if and only if its defining triple $\lbrace U,V,W \rbrace$ is
equivalent to the triple defining the algorithm of \cite{S1969}.
\end{theorem}


\begin{flushright}

\noindent Submitted to the Moscow Mathematical Society on February 23, 1972.
\end{flushright}


\end{document}